\documentclass[10pt]{iopart}
\usepackage{amsthm,amssymb,booktabs}

\theoremstyle{definition}
\newtheorem{definition}{Definition}[section]

\theoremstyle{plain}
\newtheorem{lemma}[definition]{Lemma}
\newtheorem{proposition}[definition]{Proposition}
\newtheorem{theorem}[definition]{Theorem}
\newtheorem{corollary}[definition]{Corollary}

\begin{document}

\title{The partially alternating ternary sum in an associative dialgebra}

\author{Murray R Bremner$^1$ and Juana S\'anchez Ortega$^2$}

\address{$^1$ Department of Mathematics and Statistics, University of Saskatchewan, Canada}

\address{$^2$ Department of Algebra, Geometry and Topology, University of M\'alaga, Spain}

\ead{bremner@math.usask.ca}

\ead{jsanchez@agt.cie.uma.es}



\begin{abstract}
The alternating ternary sum in an associative algebra,
  \begin{eqnarray*}
  \fl \qquad
  abc - acb - bac + bca + cab - cba,
  \end{eqnarray*}
gives rise to the partially alternating ternary sum in an associative dialgebra with products
$\dashv$ and $\vdash$ by making the argument $a$ the center of each term:
  \begin{eqnarray*}
  \fl \qquad
  a \dashv b \dashv c - a \dashv c \dashv b - b \vdash a \dashv c
  +
  c \vdash a \dashv b + b \vdash c \vdash a - c \vdash b \vdash a.
  \end{eqnarray*}
We use computer algebra to determine the polynomial identities in degree $\le 9$ satisfied by
this new trilinear operation.  In degrees 3 and 5 we obtain
  \begin{eqnarray*}
  \fl \qquad
  [a,b,c] + [a,c,b] \equiv 0,
  \quad
  [a,[b,c,d],e] + [a,[c,b,d],e] \equiv 0;
  \end{eqnarray*}
these identities define a new variety of partially alternating ternary algebras.  We show that there is
a 49-dimensional space of multilinear identities in degree 7, and we find equivalent nonlinear identities.
We use the representation theory of the symmetric group to show that there are no new identities in degree 9.
\end{abstract}


\section{Introduction}

The Lie bracket $[a,b] = ab - ba$ gives rise to two distinct ternary operations:
the Lie triple product $[[a,b],c] = abc - bac - cab + cba$,
and the alternating ternary sum (ATS) $[a,b,c] = abc - acb - bac + bca + cab - cab$,
also called the ternary commutator \cite{B} or `ternutator' \cite{DFNW}.
Apart from the obvious skew-symmetry in degree 3, the simplest non-trivial identities for the ATS have
degree 7 and were found by Bremner \cite{B} using computer algebra.  Subsequent work of Bremner and Hentzel
showed that there are no new identities for the ATS in degree 9.
(An identity is `new' if it is not a consequence of identities in lower degree.)
The ATS is closely related to the $n = 3$ case of $n$-Lie algebras introduced by Filippov \cite{F}.

The ATS appears in physics in the context of a novel formulation of quantum mechanics developed by
Nambu \cite{N}.
He introduced a multilinear $n$-bracket (now called the Nambu $n$-bracket) which becomes the ATS for
$n = 3$.
This theory has been developed by Takhtajan \cite{T}, Gautheron \cite{Ga}, Curtright and Zachos \cite{CZ},
and Ataguema et al.~\cite{AMS}, among many others.
Recently, Curtright et al.~\cite{CJM} have generalized the identity of Bremner \cite{B} to all odd $n$.
However, the most important recent developments in physics related to $n$-ary algebras are the
works of Bagger and Lambert \cite{BL} and Gustavsson \cite{Gu}, which aim at a world-volume theory of
multiple M2-branes.
For a very recent comprehensive survey of this entire area,
from both the physical and mathematical points of view,
see de Azc\'arraga and Izquierdo \cite{AI}.

Motivated by the importance of the ATS in theoretical physics,
as well as the recent works of Bremner and Peresi \cite{B-QJ, BP-SQJ}
on polynomial identities satisfied by the quasi-Jordan product in an associative dialgebra,
the present paper will focus on the partially alternating ternary sum in
an associative dialgebra.
We use computational linear algebra and the representation theory of the symmetric group
to determine the polynomial identities in degree $\le 9$ satisfied by this new trilinear operation.


\section{Preliminaries on dialgebras} \label{preliminaries}

Unless otherwise stated, the base field $F$ is the field $\mathbb{Q}$ of rational numbers.

\subsection{Dialgebras and Leibniz algebras}

Dialgebras were introduced by Loday \cite{Loday1, Loday2, Loday3} to provide a
natural setting for Leibniz algebras, a `noncommutative' generalization of Lie
algebras.

\begin{definition} (Loday \cite{Loday1}.)
A {\it Leibniz algebra} is a vector space $L$ together with a bilinear map $L \times L \to L$, denoted
$(a,b) \mapsto [a,b]$ and called the {\it Leibniz bracket}, satisfying the \emph{Leibniz identity} which
says that right multiplications are derivations:
  \begin{eqnarray*}
  \fl \qquad
  [ [ a, b ], c] \equiv [ [ a, c ], b ] + [ a, [ b, c ] ].
  \end{eqnarray*}
If $[ a, a ] \equiv 0$ then the Leibniz identity is the Jacobi identity and $L$ is a Lie algebra.
\end{definition}

Every associative algebra becomes a Lie algebra if the associative product is replaced by the Lie bracket.
Loday introduced the notion of dialgebra which gives, by a similar procedure, a Leibniz algebra: one replaces
$ab$ and $ba$ by two distinct operations, so that the resulting bracket is
not necessarily skew-symmetric.

\begin{definition}
(Loday \cite{Loday2}.)
An {\it (associative) dialgebra} is a vector space $A$ together with two bilinear maps $A \times A \to A$,
denoted $a \dashv b$ and $a \vdash b$ and called the {\it left} and {\it right} products, satisfying the
following identities:
  \begin{eqnarray*}
  \fl \qquad
  ( a \dashv b ) \dashv c \equiv a \dashv ( b \dashv c ), \quad
  ( a \vdash b ) \vdash c \equiv a \vdash ( b \vdash c ), \quad
  ( a \vdash b ) \dashv c \equiv a \vdash ( b \dashv c ),
  \\
  \fl \qquad
  ( a \dashv b ) \vdash c \equiv ( a \vdash b ) \vdash c, \quad
  a \dashv ( b \dashv c ) \equiv a \dashv ( b \vdash c ).
  \end{eqnarray*}
Both products are associative, they satisfy \emph{inner associativity} (with the operation symbols pointing
inwards), and the \emph{bar properties} (on the bar side, the operation symbols are interchangeable).
The \emph{dicommutator} is defined by $[ a, b ] = a \dashv b - b \vdash a$.
\end{definition}

Every dialgebra becomes a Leibniz algebra if the associative products are replaced by the dicommutator.

\subsection{Free dialgebras}

\begin{definition} (Loday \cite{Loday3}.)
A {\it (dialgebra) monomial} on  a set $X$ is a product $x = a_1 a_2 \cdots a_n$ where $a_1, \ldots, a_n \in X$
with some placement of parentheses and some choice of operations.
The {\it center} of $x$ is defined inductively: if $n = 1$ then $c(x) = x$; if $n \ge 2$ then
$x = y \dashv z$ or $x = y \vdash z$ and we set $c( y \dashv z ) = c(y)$ or $c( y \vdash z ) = c(z)$.
\end{definition}

A monomial is determined by the order of its factors and the position of its center.

\begin{lemma} \emph{(Loday \cite{Loday3}.)}
If $x = a_1 a_2 \cdots a_n$ is a monomial with $c(x) = a_i$ then
  \begin{eqnarray*}
  \fl \qquad
  x = ( a_1 \vdash \cdots \vdash a_{i-1} ) \vdash a_i \dashv ( a_{i+1} \dashv \cdots \dashv a_n ).
  \end{eqnarray*}
\end{lemma}

\begin{definition}
The right side of the last equation is the {\it normal form} of $x$ and
is abbreviated by the \emph{hat notation} $a_1 \cdots a_{i-1} \widehat{a}_i a_{i+1} \cdots a_n$.
\end{definition}

\begin{lemma} \emph{(Loday \cite{Loday3}.)}
The set of monomials $a_1 \cdots a_{i-1} \widehat{a}_i a_{i+1} \cdots a_n$ in normal form with
$1 \le i \le n$ and $a_1, \dots, a_n \in X$ forms a basis of the free dialgebra on $X$.
\end{lemma}

\subsection{Identities for algebras and identities for dialgebras}

Kolesnikov \cite{K} and Pozhidaev \cite{P} recently introduced an algorithm for
passing from identities of algebras to identities of dialgebras; we do not assume associativity.
Let $I$ be a multilinear algebra identity in the variables $a_1, \dots, a_n$.
For each $i = 1, \dots, n$ we convert $I$ into a multilinear dialgebra identity by making $a_i$ the center of
each monomial. In this way, one algebra identity in degree $n$ produces $n$ dialgebra identities.
For example, associativity $(ab)c - a(bc)$ gives rise to
$( \widehat{a} b ) c - \widehat{a} ( b c )$,
$( a \widehat{b} ) c - a ( \widehat{b} c )$,
$( a b ) \widehat{c} - a ( b \widehat{c} )$:
associativity of the left product, inner associativity, and associativity of the right product.

The same algorithm can be used to convert a multilinear algebra operation into a set of multilinear
dialgebra operations.
For example, the Lie bracket $ab - ba$ gives rise to the left and right Leibniz products
$a \dashv b - b \vdash a$ and $a \vdash b - b \dashv a$.
We use this method to obtain ternary dialgebra operations from ternary algebra operations.


\section{Degree 3}

The Kolesnikov-Pozhidaev algorithm produces the dialgebra version of the ATS.

\begin{definition} \label{opdef}
The \textit{partially alternating ternary sum} (\textit{PATS}) is this trilinear operation
in an associative dialgebra:
  \begin{eqnarray*}
  \fl \qquad
  [a,b,c]
  =
  \widehat{a}bc
  - \widehat{a}cb
  - b\widehat{a}c
  + c\widehat{a}b
  + bc\widehat{a}
  - cb\widehat{a}.
  \end{eqnarray*}
It is easy to see that the PATS is skew-symmetric in its second and third arguments;
Proposition \ref{degree5} shows that it becomes completely alternating in the second and third arguments of a nested monomial.
The PATS is obtained from the ATS by making $a$ the center of each monomial;
the skew-symmetry of the ATS implies that making $b$ or $c$ the center gives equivalent operations.
\end{definition}

\begin{definition}  \label{Pdefinition}
A \emph{ternary algebra} is a vector space $T$ together with a trilinear map $T \times T \times T \to T$
denoted by $(a, b, c)$.  In a ternary algebra we define the polynomial $P(a,b,c)= (a,b,c) + (a,c,b)$.
\end{definition}

\begin{proposition} \label{degree3}
Every multilinear polynomial identity in degree $\le 3$ satisfied by the PATS is a consequence
of $P(a,b,c) \equiv 0$.
\end{proposition}

\begin{proof}
Consider the general trilinear identity of degree 3:
  \begin{eqnarray*}
  \fl \qquad
    x_1 [a,b,c]
  + x_2 [a,c,b]
  + x_3 [b,a,c]
  + x_4 [b,c,a]
  + x_5 [c,a,b]
  + x_6 [c,b,a]
  =
  0.
  \end{eqnarray*}
Each of the 6 ternary monomials expands using the PATS into a linear combination of 6 dialgebra monomials.
Altogether we obtain 18 dialgebra monomials:
  \begin{eqnarray*}
  \fl \qquad
  \widehat{a}bc, \widehat{a}cb, \widehat{b}ac,
  \widehat{b}ca, \widehat{c}ab, \widehat{c}ba, \,
  a\widehat{b}c, a\widehat{c}b, b\widehat{a}c,
  b\widehat{c}a, c\widehat{a}b, c\widehat{b}a, \,
  ab\widehat{c}, ac\widehat{b}, ba\widehat{c},
  bc\widehat{a}, ca\widehat{b}, cb\widehat{a}.
  \end{eqnarray*}
Let $E$ be the $18 \times 6$ expansion matrix whose $(i,j)$ entry is the coefficient of the $i$-th
dialgebra monomial in the expansion of the $j$-th ternary monomial (Table \ref{expansion3}).
The coefficient vectors of the identities in degree 3 satisfied by the PATS are the vectors in the nullspace
of $E$. We compute the row canonical form of $E$ and find that the canonical integral basis of the nullspace
consists of three permutations of $P(a,b,c)$, namely $[a,b,c] + [a,c,b]$, $[b,a,c] + [b,c,a]$ and
$[c,a,b] + [c,b,a]$.
\end{proof}

  \begin{table}
  \caption{Transpose of the expansion matrix $E$ in degree 3.}
  \label{expansion3}
  \begin{indented}
  \item[]
  \begin{tabular}{c}
  \br
  $
  \left[
  \begin{array}{rrrrrrrrrrrrrrrrrr}
   1 &\!\!\!\! -1 &\!\!\!\!  . &\!\!\!\!  . &\!\!\!\!  . &\!\!\!\!  . &\!\!\!\!
   . &\!\!\!\!  . &\!\!\!\! -1 &\!\!\!\!  . &\!\!\!\!  1 &\!\!\!\!  . &\!\!\!\!
   . &\!\!\!\!  . &\!\!\!\!  . &\!\!\!\!  1 &\!\!\!\!  . &\!\!\!\! -1
  \\
  -1 &\!\!\!\!  1 &\!\!\!\!  . &\!\!\!\!  . &\!\!\!\!  . &\!\!\!\!  . &\!\!\!\!
   . &\!\!\!\!  . &\!\!\!\!  1 &\!\!\!\!  . &\!\!\!\! -1 &\!\!\!\!  . &\!\!\!\!
   . &\!\!\!\!  . &\!\!\!\!  . &\!\!\!\! -1 &\!\!\!\!  . &\!\!\!\!  1
  \\
   . &\!\!\!\!  . &\!\!\!\!  1 &\!\!\!\! -1 &\!\!\!\!  . &\!\!\!\!  . &\!\!\!\!
  -1 &\!\!\!\!  . &\!\!\!\!  . &\!\!\!\!  . &\!\!\!\!  . &\!\!\!\!  1 &\!\!\!\!
   . &\!\!\!\!  1 &\!\!\!\!  . &\!\!\!\!  . &\!\!\!\! -1 &\!\!\!\!  .
  \\
   . &\!\!\!\!  . &\!\!\!\! -1 &\!\!\!\!  1 &\!\!\!\!  . &\!\!\!\!  . &\!\!\!\!
   1 &\!\!\!\!  . &\!\!\!\!  . &\!\!\!\!  . &\!\!\!\!  . &\!\!\!\! -1 &\!\!\!\!
   . &\!\!\!\! -1 &\!\!\!\!  . &\!\!\!\!  . &\!\!\!\!  1 &\!\!\!\!  .
  \\
   . &\!\!\!\!  . &\!\!\!\!  . &\!\!\!\!  . &\!\!\!\!  1 &\!\!\!\! -1 &\!\!\!\!
   . &\!\!\!\! -1 &\!\!\!\!  . &\!\!\!\!  1 &\!\!\!\!  . &\!\!\!\!  . &\!\!\!\!
   1 &\!\!\!\!  . &\!\!\!\! -1 &\!\!\!\!  . &\!\!\!\!  . &\!\!\!\!  .
  \\
   . &\!\!\!\!  . &\!\!\!\!  . &\!\!\!\!  . &\!\!\!\! -1 &\!\!\!\!  1 &\!\!\!\!
   . &\!\!\!\!  1 &\!\!\!\!  . &\!\!\!\! -1 &\!\!\!\!  . &\!\!\!\!  . &\!\!\!\!
  -1 &\!\!\!\!  . &\!\!\!\!  1 &\!\!\!\!  . &\!\!\!\!  . &\!\!\!\!  .
  \end{array}
  \right]
  $
  \end{tabular}

  \end{indented}
  \end{table}


\section{Degree 5}

For a ternary operation, there are three association types in degree 5:
$((a,b,c),d,e)$, $(a,(b,c,d),e)$ and $(a,b,(c,d,e))$.

\begin{lemma} \label{degree5monomials}
If a ternary operation satisfies $P(a,b,c) \equiv 0$, then every multilinear monomial in degree 5 equals
one of the following 90 monomials:
  \begin{eqnarray*}
  \fl \qquad
  ((a^\sigma,b^\sigma,c^\sigma),d^\sigma,e^\sigma)
  \;\; \mathrm{for} \;\;
  b^\sigma < c^\sigma, \, d^\sigma < e^\sigma,
  \quad
  (a^\sigma,(b^\sigma,c^\sigma,d^\sigma),e^\sigma)
  \;\; \mathrm{for} \;\;
  c^\sigma < d^\sigma.
  \end{eqnarray*}
Here $\sigma$ is a permutation of $\{a,b,c,d,e\}$ and $<$ denotes alphabetical precedence.
\end{lemma}

\begin{proof}
Since $P(a,b,c) \equiv 0$ implies $(a,b,(c,d,e)) = -(a,(c,d,e),b)$, we can ignore the third association type.
Applying $P(a,b,c) \equiv 0$ to the first and second types gives
  \begin{eqnarray*}
  \fl \qquad
  ((a,b,c),d,e)
  &
  \equiv -((a,c,b),d,e) \equiv -((a,b,c),e,d) \equiv ((a,c,b),e,d),
  \\
  \fl \qquad
  (a,(b,c,d),e)
  &
  \equiv -(a,(b,d,c),e).
  \end{eqnarray*}
Hence, there are $5!/4 = 30$ monomials in the first type and $5!/2 = 60$ in the second.
\end{proof}

\begin{definition} \label{Qdefinition}
In a ternary algebra, we define the polynomial $Q(a,b,c,d,e) = (a,(b,c,d),e) + (a,(c,b,d),e)$.
\end{definition}

\begin{proposition} \label{degree5}
Every multilinear polynomial identity in degree $\le 5$ satisfied by the PATS
is a consequence of $P(a,b,c) \equiv 0$ and $Q(a,b,c,d,e) \equiv 0$.
\end{proposition}

\begin{proof}
We order the multilinear ternary monomials of Lemma \ref{degree5monomials} first by association type
and then by lex order of the permutation. Each ternary monomial expands into a sum of 36 dialgebra monomials.
For the first two association types we have:
  \begin{eqnarray*}
  \fl \qquad
  [[a,b,c],d,e] =
  \\
  \fl \qquad
  \widehat{a} b c d e
  - \widehat{a} b c e d
  - \widehat{a} c b d e
  + \widehat{a} c b e d
  - b \widehat{a} c d e
  + b \widehat{a} c e d
  + c \widehat{a} b d e
  - c \widehat{a} b e d
  - d \widehat{a} b c e
  \\
  \fl \qquad
  {}
  + d \widehat{a} c b e
  + e \widehat{a} b c d
  - e \widehat{a} c b d
  + b c \widehat{a} d e
  - b c \widehat{a} e d
  - c b \widehat{a} d e
  + c b \widehat{a} e d
  + d b \widehat{a} c e
  - d c \widehat{a} b e
  \\
  \fl \qquad
  {}
  + d e \widehat{a} b c
  - d e \widehat{a} c b
  - e b \widehat{a} c d
  + e c \widehat{a} b d
  - e d \widehat{a} b c
  + e d \widehat{a} c b
  - d b c \widehat{a} e
  + d c b \widehat{a} e
  - d e b \widehat{a} c
  \\
  \fl \qquad
  {}
  + d e c \widehat{a} b
  + e b c \widehat{a} d
  - e c b \widehat{a} d
  + e d b \widehat{a} c
  - e d c \widehat{a} b
  + d e b c \widehat{a}
  - d e c b \widehat{a}
  - e d b c \widehat{a}
  + e d c b \widehat{a},
  \\
  \fl \qquad
  {}
  [a,[b,c,d],e] =
  \\
  \fl \qquad
    \widehat{a} b c d e
  - \widehat{a} b d c e
  - \widehat{a} c b d e
  + \widehat{a} c d b e
  + \widehat{a} d b c e
  - \widehat{a} d c b e
  - \widehat{a} e b c d
  + \widehat{a} e b d c
  + \widehat{a} e c b d
  \\
  \fl \qquad
  {}
  - \widehat{a} e c d b
  - \widehat{a} e d b c
  + \widehat{a} e d c b
  + e \widehat{a} b c d
  - e \widehat{a} b d c
  - e \widehat{a} c b d
  + e \widehat{a} c d b
  + e \widehat{a} d b c
  - e \widehat{a} d c b
  \\
  \fl \qquad
  {}
  - b c d \widehat{a} e
  + b d c \widehat{a} e
  + c b d \widehat{a} e
  - c d b \widehat{a} e
  - d b c \widehat{a} e
  + d c b \widehat{a} e
  + b c d e \widehat{a}
  - b d c e \widehat{a}
  - c b d e \widehat{a}
  \\
  \fl \qquad
  {}
  + c d b e \widehat{a}
  + d b c e \widehat{a}
  - d c b e \widehat{a}
  - e b c d \widehat{a}
  + e b d c \widehat{a}
  + e c b d \widehat{a}
  - e c d b \widehat{a}
  - e d b c \widehat{a}
  + e d c b \widehat{a}.
  \end{eqnarray*}
Altogether there are $5 \cdot 5! = 600$ dialgebra monomials; we order them first by the position of the
center and then
lexicographically. We construct the expansion matrix $E$ of size $600 \times 90$ whose $(i,j)$ entry is
the coefficient of the $i$-th dialgebra monomial in the expansion of the $j$-th ternary monomial.
We use the \texttt{Maple} package \texttt{LinearAlgebra} to manipulate this matrix: \texttt{Rank} returns
the value 50 and so the nullspace has dimension 40; \texttt{ReducedRowEchelonForm} computes the row canonical
form; \texttt{DeleteRow} removes the 550 zero rows at the bottom of the matrix.
The result is a $50 \times 90$ matrix with the leading 1s of its rows in columns
1--33, 36, 43--45, 48, 55--57, 60, 67--69, 72, 79--81, 84,
and so the free variables correspond to columns
34, 35, 37--42, 46, 47, 49--54, 58, 59, 61--66, 70, 71, 73--78, 82, 83, 85--90.
For each free variable, we set that variable to 1 and the other free variables to 0, and solve for
the leading variables.  We obtain a basis of the nullspace consisting of the following identities,
20 of each form:
  \begin{eqnarray*}
  \fl \qquad
  \left.
  \begin{array}{l}
  [a^\sigma,[b^\sigma,c^\sigma,d^\sigma],e^\sigma]
  +
  [a^\sigma,[c^\sigma,b^\sigma,d^\sigma],e^\sigma]
  \equiv 0
  \\
  {}
  [a^\sigma,[b^\sigma,c^\sigma,d^\sigma],e^\sigma ]
  -
  [a^\sigma,[d^\sigma,b^\sigma,c^\sigma],e^\sigma]
  \equiv 0
  \end{array}
  \right\}
  \; \mathrm{for} \;
  b^\sigma < c^\sigma < d^\sigma.
  \end{eqnarray*}
Every identity of the first form is equivalent to $Q(a,b,c,d,e) \equiv 0$. We have
  \begin{eqnarray*}
  \fl \qquad
  [a^\sigma,[b^\sigma,c^\sigma,d^\sigma],e^\sigma]
  \stackrel{P}{\equiv}
  -[a^\sigma,[b^\sigma,d^\sigma,c^\sigma],e^\sigma]
  \stackrel{Q}{\equiv}
  [a^\sigma,[d^\sigma,b^\sigma,c^\sigma],e^\sigma],
  \end{eqnarray*}
and so every identity of the second form follows from $P$ and $Q$.
\end{proof}

\begin{corollary}
If a ternary algebra satisfies $P(a,b,c) \equiv 0$ and $Q(a,b,c,d,e) \equiv 0$,
then every monomial $(a,(b,c,d),e)$ is an alternating function of $b, c, d$
and every monomial $(a,b,(c,d,e))$ is an alternating function of $c,d,e$.
\end{corollary}


\section{Partially alternating ternary algebras} \label{PAsection}

The PATS makes any dialgebra into a partially alternating ternary algebra as follows.

\begin{definition}
A \emph{completely alternating ternary algebra (CATA)} is one with a \emph{completely alternating product}:
$(a^\sigma, b^\sigma, c^\sigma) = \epsilon(\sigma) (a,b,c)$ for all permutations $\sigma$ where $\epsilon$
is the sign.
A \emph{partially alternating ternary algebra (PATA)} is one with a
\emph{partially alternating product}: the product satisfies
$P(a,b,c) \equiv 0$ and $Q(a,b,c,d,e) \equiv 0$.
\end{definition}

\begin{lemma}
Let $(\cdot, \cdot, \cdot)$ be a partially alternating ternary product and let $s = a_1 \cdots a_n$ be a
monomial with some placement of (ternary) parentheses.
Then the second and third arguments of any submonomial $(x,y,z)$ of $s$ are completely alternating.
\end{lemma}

\begin{proof}
We proceed by induction on $n$. For $n \le 3$ the claim is obvious.
Set $n > 3$ and assume that the result holds for any monomial of degree $< n$. We have the factorization
$s = (t,u,v)$; then the second or third argument of a submonomial of $s$ is
either ($i$) $u$ or $v$, or ($ii$) the second or third argument of a submonomial of $t$, $u$ or $v$.
In case ($ii$), the claim follows from the inductive hypothesis.
In case ($i$), it suffices to prove the claim for the second argument $u$,
since $P \equiv 0$ implies $s = -(t,v,u)$ and so the claim also holds for $v$.
To show the claim for $u$, it is enough to note that $P \equiv 0$ (respectively $Q \equiv 0$)
implies that $u$ alternates in its second and third arguments (respectively its first and second arguments).
Since these two transpositions generate the symmetric group, $u$ is completely alternating.
\end{proof}

A ternary operation has monomials only in odd degrees. We generate the association types for a partially
alternating ternary product inductively by degree; we simultaneously generate the completely alternating
(CA) types and the partially alternating (PA) types (Table \ref{CAPAtypes}). Suppose that we have already
generated ordered lists of the CA and PA types in degrees $< n$. To generate the CA and PA association
types of degree $n$, we consider all partitions of $n = i + j + k$ into three odd parts:
  \begin{enumerate}
  \item
  If $i > j > k$ then for all CA types $t, u, v$ of degrees $i, j, k$ respectively, we include $(t,u,v)$ as
  a new CA type in degree $n$.
  \item
  If $i = j > k$ then for all CA types $t, u, v$ of degrees $i, i, k$ where $t$ precedes $u$ in the types of
  degree $i$, we include $(t,u,v)$ as a CA type in degree $n$.
  \item
  If $i > j = k$ then for all CA types $t, u, v$ of degrees $i, j, j$ where $u$ precedes $v$ in the types of
  degree $j$, we include $(t,u,v)$ as a CA type in degree $n$.
  \item
  If $i = j = k$ then for all CA types $t, u, v$ of degree $i$ where $t$ precedes $u$ and $u$ precedes $v$,
  we include $(t,u,v)$ as a CA type in degree $n$.
  \item
  If $j > k$ then for all PA types $t$ of degree $i$ and all CA types $u, v$ of degrees $j, k$ respectively,
  we include $(t,u,v)$ as a new PA type in degree $n$.
  \item
  If $j = k$ then for all PA types $t$ of degree $i$ and all CA types $u, v$ of degree $j$ where $u$ precedes
  $v$, we include $(t,u,v)$ as a PA type in degree $n$.
  \end{enumerate}
To enumerate the multilinear monomials in a completely or partially alternating association type in degree
$n$, we need to determine the number $d$ of skew-symmetries of the association type; the number of monomials
is then $n!/d$. To compute $d$ we use the recursive procedure of Figure \ref{countsymmetry}; before calling
the procedure we set $d \leftarrow 1$.

  \begin{table}
  \caption{Association types for ternary products.}
  \label{CAPAtypes}
  \begin{indented}
  \item[]
  \begin{tabular}{@{}ccccc}
  \br
  $a$ & $(a,b,c)$ & $((a,b,c),d,e)$ & $(((a,b,c),d,e),f,g)$ & $((((a,b,c),d,e),f,g),h,i)$ \\
      &           &                 & $((a,b,c),(d,e,f),g)$ & $(((a,b,c),(d,e,f),g),h,i)$ \\
      &           &                 &                       & $(((a,b,c),d,e),(f,g,h),i)$ \\
  \multicolumn{4}{l}{completely alternating types}          & $((a,b,c),(d,e,f),(g,h,i))$ \\
  \mr
  $a$ & $(a,b,c)$ & $((a,b,c),d,e)$ & $(((a,b,c),d,e),f,g)$ & $((((a,b,c),d,e),f,g),h,i)$ \\
      &           & $(a,(b,c,d),e)$ & $((a,(b,c,d),e),f,g)$ & $(((a,(b,c,d),e),f,g),h,i)$ \\
      &           &                 & $((a,b,c),(d,e,f),g)$ & $(((a,b,c),(d,e,f),g),h,i)$ \\
      &           &                 & $(a,((b,c,d),e,f),g)$ & $((a,((b,c,d),e,f),g),h,i)$ \\
      &           &                 & $(a,(b,c,d),(e,f,g))$ & $((a,(b,c,d),(e,f,g)),h,i)$ \\
      &           &                 &                       & $(((a,b,c),d,e),(f,g,h),i)$ \\
      &           &                 &                       & $((a,(b,c,d),e),(f,g,h),i)$ \\
      &           &                 &                       & $((a,b,c),((d,e,f),g,h),i)$ \\
      &           &                 &                       & $((a,b,c),(d,e,f),(g,h,i))$ \\
      &           &                 &                       & $(a,(((b,c,d),e,f),g,h),i)$ \\
      &           &                 &                       & $(a,((b,c,d),(e,f,g),h),i)$ \\
  \multicolumn{4}{l}{partially alternating types}           & $(a,((b,c,d),e,f),(g,h,i))$ \\
  \br
  \end{tabular}
  \end{indented}
  \end{table}

  \begin{figure}
  \caption{Recursive procedure $\texttt{countsymmetry}(x,\texttt{flag})$.}
  \label{countsymmetry}
  \begin{indented}
  \item[]
  \begin{tabular}{l}
  \br
  \emph{Input:}
  \\
  A completely or partially alternating ternary association type $x$;
  \\
  a Boolean variable \texttt{flag} which is \texttt{true} for CA and \texttt{false} for PA.
  \\
  \mr
  \emph{Procedure:}
  \\
  If $\deg(x) > 1$ then write $x = (x_1,x_2,x_3)$:
  \\
  \quad If $\texttt{flag} = \texttt{true}$ then
  \\
  \quad \quad \texttt{countsymmetry}( $x_1$, \texttt{true} ):
  \\
  \quad \quad \texttt{countsymmetry}( $x_2$, \texttt{true} ):
  \\
  \quad \quad \texttt{countsymmetry}( $x_3$, \texttt{true} ):
  \\
  \quad \quad If all three of $x_1, x_2, x_3$ have the same degree and association type then
  \\
  \quad \quad \quad
  set $d \leftarrow 6d$.
  \\
  \quad \quad If only two of $x_1, x_2, x_3$ have the same degree and association type then
  \\
  \quad \quad \quad set $d \leftarrow 2d$.
  \\
  \quad else
  \\
  \quad \quad \texttt{countsymmetry}( $x_1$, \texttt{false} ):
  \\
  \quad \quad \texttt{countsymmetry}( $x_2$, \texttt{true} ):
  \\
  \quad \quad \texttt{countsymmetry}( $x_3$, \texttt{true} ):
  \\
  \quad \quad If $x_2$, $x_3$ have the same degree and association type then set $d \leftarrow 2d$.
  \\
  \br
  \end{tabular}
  \end{indented}
  \end{figure}


\section{Degree 7} \label{degree7section}

From now on we usually omit the commas in all ternary monomials.

\begin{lemma} \label{degree7monomials}
In a partially alternating ternary algebra, every multilinear monomial in degree 7 equals one of the
following 1960 monomials:
  \begin{eqnarray*}
  \fl \qquad
  (1)
  &\qquad
  (((a^\sigma b^\sigma c^\sigma) d^\sigma e^\sigma) f^\sigma g^\sigma)
  &\qquad
  (b^\sigma < c^\sigma, d^\sigma < e^\sigma, f^\sigma < g^\sigma)
  \\
  \fl \qquad
  (2)
  &\qquad
  ((a^\sigma (b^\sigma c^\sigma d^\sigma) e^\sigma) f^\sigma g^\sigma)
  &\qquad
  (b^\sigma < c^\sigma < d^\sigma, f^\sigma < g^\sigma)
  \\
  \fl \qquad
  (3)
  &\qquad
  ((a^\sigma b^\sigma c^\sigma) (d^\sigma e^\sigma f^\sigma) g^\sigma)
  &\qquad
  (b^\sigma < c^\sigma, d^\sigma < e^\sigma < f^\sigma)
  \\
  \fl \qquad
  (4)
  &\qquad
  (a^\sigma ((b^\sigma c^\sigma d^\sigma) e^\sigma f^\sigma) g^\sigma)
  &\qquad
  (b^\sigma < c^\sigma < d^\sigma, e^\sigma < f^\sigma)
  \\
  \fl \qquad
  (5)
  &\qquad
  (a^\sigma (b^\sigma c^\sigma d^\sigma) (e^\sigma f^\sigma g^\sigma))
  &\qquad
  (b^\sigma < c^\sigma < d^\sigma, e^\sigma < f^\sigma < g^\sigma, b^\sigma < e^\sigma)
  \end{eqnarray*}
Here $\sigma$ is a permutation of $\{a,b,c,d,e,f,g\}$ and $<$ is alphabetical precedence.
\end{lemma}

\begin{proof}
Degree 7 has 12 ternary association types, but in a PATA these reduce to 5:
  \begin{eqnarray*}
  \fl \qquad
  \begin{array}{ll}
  (((abc)de)fg) = \mathrm{type \; (1)},
  &\qquad
  ((a(bcd)e)fg) = \mathrm{type \; (2)},
  \\
  ((ab(cde))fg) = -((a(cde)b)fg),
  &\qquad
  ((abc)(def)g) = \mathrm{type \; (3)},
  \\
  ((abc)d(efg)) = -((abc)(efg)d),
  &\qquad
  (a((bcd)ef)g) = \mathrm{type \; (4)},
  \\
  (a(b(cde)f)g) = -(a((cde)bf)g),
  &\qquad
  (a(bc(def))g) = (a((def)bc)g),
  \\
  (a(bcd)(efg)) = \mathrm{type \; (5)}.
  &\qquad
  (ab((cde)fg)) = -(a((cde)fg)b),
  \\
  (ab(c(def)g)) = (a((def)cg)b),
  &\qquad
  (ab(cd(efg))) = -(a((efg)cd)b).
  \end{array}
  \end{eqnarray*}
To enumerate the multilinear monomials in each type, we count the skew-symmetries:
  \begin{eqnarray*}
  \fl \qquad
  \begin{tabular}{llr}
  $(((abc)de)fg)$
  &\;
  alternates in $b,c$ and $d,e$ and $f,g$
  &\;
  $7!/8 = 630$
  \\
  $((a(bcd)e)fg)$
  &\;
  alternates in $b,c,d$ and $f,g$
  &\;
  $7!/12 = 420$
  \\
  $((abc)(def)g)$
  &\;
  alternates in $b,c$ and $d,e,f$
  &\;
  $7!/12 = 420$
  \\
  $(a((bcd)ef)g)$
  &\;
  alternates in $b,c,d$ and $e,f$
  &\;
  $7!/12 = 420$
  \\
  $(a(bcd)(efg))$
  &\;
  alternates in $b,c,d$ and $e,f,g$ and $bcd, efg$
  &\;
  $7!/72 = 70$
  \end{tabular}
  \end{eqnarray*}
We order these monomials first by association type and then lexicographically.
\end{proof}

\begin{definition} \label{RSdefinition}
In a PATA we consider the following polynomials of degree 7:
  \begin{eqnarray*}
  \fl \qquad
  R(a,b,c,d,e,f,g) =
  \\
  \fl \qquad
  \quad
  \frac1{12} \sum_{\sigma \in S_6}
  \epsilon(\sigma) \,
  ((a (b^\sigma c^\sigma d^\sigma) e^\sigma) f^\sigma g^\sigma)
  -
  \frac1{12} \sum_{\sigma \in S_6}
  \epsilon(\sigma) \,
  ((a b^\sigma c^\sigma) (d^\sigma e^\sigma f^\sigma) g^\sigma),
  \\
  \fl \qquad
  S(a,b,c,d,e,f,g) =
  \\
  \fl \qquad
  \quad
  \frac14
  \sum_{\sigma \in S_5}
  \epsilon(\sigma) \,
  (((a b^\sigma c^\sigma) d^\sigma g) e^\sigma f^\sigma)
  -
  \frac16
  \sum_{\sigma \in S_5}
  \epsilon(\sigma) \,
  ( ( a ( b^\sigma c^\sigma d^\sigma ) e^\sigma ) f^\sigma g )
  \\
  \fl \qquad
  \quad
  +
  \frac14
  \sum_{\sigma \in S_5}
  \epsilon(\sigma) \,
  ( ( a ( b^\sigma c^\sigma g ) d^\sigma ) e^\sigma f^\sigma )
  +
  \frac{1}{12}
  \sum_{\sigma \in S_5}
  \epsilon(\sigma) \,
  ( ( a b^\sigma c^\sigma ) ( d^\sigma e^\sigma f^\sigma ) g )
  \\
  \fl \qquad
  \quad
  -
  \frac16
  \sum_{\sigma \in S_5}
  \epsilon(\sigma) \,
  ( a (( b^\sigma c^\sigma d^\sigma ) e^\sigma g ) f^\sigma )
  -
  \frac{1}{12}
  \sum_{\sigma \in S_5}
  \epsilon(\sigma) \,
  ( a ( b c^\sigma d^\sigma ) ( e^\sigma f^\sigma g^\sigma) ).
  \end{eqnarray*}
Some of the permutations $\sigma$ in $R$ and $S$ produce monomials which are not `straightened': $P$ and $Q$
need to be applied to convert such a monomial into ($\pm$) an equivalent monomial in which the permutation
$\sigma'$ satisfies $\sigma' < \sigma$ in lex order.  This produces repetitions; to cancel the resulting
coefficients, we use appropriate fractions. (This `straightening' algorithm is described in detail in
Subsection \ref{subsectionstraightening}.) After this process, each monomial has coefficient $\pm 1$, and
each identity has 120 terms.
\end{definition}

\begin{theorem} \label{degree7theorem}
Every multilinear polynomial identity of degree $\le 7$ satisfied by the PATS is a consequence of
$P \equiv 0$, $Q \equiv 0$, $R \equiv 0$ and $S \equiv 0$.
\end{theorem}

\begin{proof}
There are $7 \cdot 7! = 35280$ dialgebra monomials in degree 7; this and Lemma \ref{degree7monomials} show
that the expansion matrix $E$ has size $35280 \times 1960$. It is not practical to use rational arithmetic
with such a large matrix so we use \texttt{LinearAlgebra[Modular]}. The vector space of multilinear ternary
polynomials is a representation of the symmetric group which acts by permuting the variables. The group
algebra $F S_n$ is semisimple over a field $F$ of characteristic 0 or $p > n$, and so for degree 7 any $p > 7$
will give the correct ranks. The procedure \texttt{RowReduce} computes the row canonical form of $E$ and its
rank 1911. The 49-dimensional nullspace consists of the multilinear polynomial identities in degree 7 for the
PATS which are not consequences of $P \equiv 0$ and $Q \equiv 0$. The procedure \texttt{Basis} computes the
canonical basis of the nullspace, which we sort by increasing number of terms in the corresponding identities:
  \begin{eqnarray*}
  \fl \qquad
  \begin{tabular}{l}
  The first 7 identities: 120 terms; coefficients $\pm 1$; association types 2, 3.
  \\
  The next 28: 120 terms; coefficients $\pm 1$; association types 1, 2, 3, 4.
  \\
  The next 7: 120 terms; coefficients $\pm 1$; association types 1, 2, 3, 4, 5.
  \\
  The last 7: 180 terms; coefficients $\pm 1$, $\pm 2$; association types 1, 2, 3, 4.
  \end{tabular}
  \end{eqnarray*}
Further computations verify that identities 1 and 36 ($R$ and $S$) generate the nullspace: every identity is
a linear combination of permutations of $R$ and $S$. Moreover, $R$ (respectively, $S$) generates a subspace
of dimension 7 (respectively, 42). Hence the nullspace is the direct sum of these two subspaces, and so $R$
and $S$ are independent.
\end{proof}


\section{Representation theory of the symmetric group}

To study nonlinear identities, and identities of higher degree, we use the representation theory of
the symmetric group $S_n$; our main reference is James and Kerber \cite{JK}.

\begin{definition}
A \emph{partition} of $n$ is a tuple $\lambda = (n_1,\dots,n_\ell)$ with $n = n_1 + \cdots + n_\ell$ and
$n_1 \ge \cdots \ge n_\ell \ge 1$. The \emph{frame} $[\lambda]$ consists of $n$ boxes in $\ell$ left-justified
rows with $n_i$ boxes in row $i$. A \emph{tableau} is a bijection between $\{ 1, \dots, n \}$ and the boxes of
$[\lambda]$. In a \emph{standard} tableau the numbers increase from left to right and from top to bottom.
\end{definition}

The irreducible representations of $S_n$ correspond to the partitions of $n$. The dimension $d_\lambda$
associated to $\lambda$ is the number of standard tableaux
with frame $[\lambda]$. If $F = \mathbb{Q}$ or $F = \mathbb{F}_p$ for $p > n$ then the group algebra
$F S_n$ decomposes into an orthogonal direct sum of two-sided ideals
isomorphic to simple matrix algebras:
  \begin{eqnarray*}
  \fl \qquad
  F S_n \approx \bigoplus_{\lambda} \, M_{d_\lambda}(F).
  \end{eqnarray*}
Given a permutation $\pi$ and a partition $\lambda$, we need to compute the projection of $\pi$ onto
$M_{d_\lambda}(F)$: the matrix for $\pi$ in representation $\lambda$. A simple algorithm for this was found
by Clifton \cite{C}. We fix an ordering of the standard tableaux $T_1, \dots, T_d$ ($d = d_\lambda$) and
construct a matrix $R^\lambda_\pi$ as follows:  Apply $\pi$ to $T_j$, obtaining a (possibly non-standard)
tableau $\pi T_j$. If there exist two numbers in the same column of $T_i$ and the same row of $\pi T_j$,
then $(R^\lambda_\pi)_{ij} = 0$. Otherwise, $(R^\lambda_\pi)_{ij}$ is the sign of the permutation of
$T_i$ which leaves the columns invariant as sets and moves the numbers to the correct rows of $\pi T_j$.
The matrix $R^\lambda_\mathrm{id}$ for the identity permutation may not be the identity matrix but is
invertible. An algorithm for computing $R^\lambda_\pi$ is given by Bremner and Peresi \cite{BP-SQJ}.

\begin{lemma}
\emph{(Clifton \cite{C}.)}
\label{cliftonlemma}
The matrix representing $\pi$ in partition $\lambda$ equals $(R^\lambda_\mathrm{id})^{-1} R^\lambda_\pi$.
\end{lemma}

Any polynomial identity (not necessarily multilinear or even homogeneous) of degree $\le n$ over a field $F$
of characteristic 0 or $p > n$ is equivalent to a finite set of multilinear identities \cite{Z}. We consider
a multilinear identity $I( x_1, \dots, x_n )$ of degree $n$ and collect the terms with the same association
type: $I = I_1 + \cdots + I_t$. The monomials in each $I_k$ differ only by a permutation of $x_1, \dots, x_n$;
hence $I_k$ is an element of $F S_n$ and $I$ is an element of $( F S_n )^t$. If $U\subseteq ( F S_n )^t$ is
the span of all permutations of $I$, then $U$ is a representation of $S_n$, and so $U$ is the direct sum of
components corresponding to the irreducible representations of $S_n$.  This breaks down a large computational
problem into smaller pieces. We fix a partition $\lambda$ with associated dimension $d = d_\lambda$. To
determine the $\lambda$-component of $U$, we construct a $d \times dt$ matrix $M_\lambda$ consisting of $t$
blocks of size $d \times d$; in block $j$ we put the representation matrix for the terms of $I_j$.

\begin{definition}
The rank of $M_\lambda$ is the \emph{rank of identity $I$ in partition $\lambda$}.
\end{definition}

  \begin{table}
  \caption{Representation matrix $X_\lambda$ in partition $\lambda$.}
  \label{repmatexp}
  \begin{indented}
  \item[]
  \begin{tabular}{@{}l}
  \br
  $
  \left[
  \begin{array}{ccccc|ccccc}
  \rho_\lambda(E^1_1) &\!\!\!\! \rho_\lambda(E^1_2) &\!\!\!\! \cdots &\!\!\!\!
  \rho_\lambda(E^1_{n-1}) &\!\!\!\! \rho_\lambda(E^1_n) &\!
  -I_d   &\!\!\!\!  O     &\!\!\!\! \cdots &\!\!\!\!  O     &\!\!\!\!  O
  \\
  \rho_\lambda(E^2_1) &\!\!\!\! \rho_\lambda(E^2_2) &\!\!\!\! \cdots &\!\!\!\!
  \rho_\lambda(E^2_{n-1}) &\!\!\!\! \rho_\lambda(E^2_n) &\!
  O      &\!\!\!\! -I_d   &\!\!\!\! \cdots &\!\!\!\!  O     &\!\!\!\!  O
  \\
  \vdots &\!\!\!\! \vdots &\!\!\!\! \ddots &\!\!\!\! \vdots &\!\!\!\! \vdots &\!
  \vdots &\!\!\!\! \vdots &\!\!\!\! \ddots &\!\!\!\! \vdots &\!\!\!\! \vdots
  \\
  \rho_\lambda(E^{t-1}_1) &\!\!\!\! \rho_\lambda(E^{t-1}_2) &\!\!\!\! \cdots &\!\!\!\!
  \rho_\lambda(E^{t-1}_{n-1}) &\!\!\!\! \rho_\lambda(E^{t-1}_n) &\!
  O      &\!\!\!\! O      &\!\!\!\! \cdots &\!\!\!\! -I_d   &\!\!\!\!  O
  \\
  \rho_\lambda(E^t_1) &\!\!\!\! \rho_\lambda(E^t_2) &\!\!\!\! \cdots &\!\!\!\!
  \rho_\lambda(E^t_{n-1}) &\!\!\!\! \rho_\lambda(E^t_n) &\!
  O      &\!\!\!\! O      &\!\!\!\! \cdots &\!\!\!\!  O     &\!\!\!\! -I_d
  \end{array}
  \right]
  $
  \end{tabular}
  \end{indented}
  \end{table}

We modify this procedure to determine the nullspace of the expansion matrix $E$ for the PATS. In degree $n$,
let $t = t_n$ be the number of association types for a partially alternating ternary product. Consider the
monomial in ternary association type $i$ with the identity permutation of the variables, and let $E^i$ be its
expansion using the PATS. We have $E^i = E^i_1 + \cdots + E^i_n$ where $E^i_j$ contains the dialgebra monomials
with center in position $j$.  We construct a $\, td \times (n{+}t)d\, $ matrix $X_\lambda$ with $t$ rows and
$n{+}t$ columns of $d \times d$ blocks (Table \ref{repmatexp}). In the right side, in block $(i,n{+}i)$ for
$1 \le i \le t$, we put $-I_d$ (identity matrix); the other blocks of the right side are zero. In the left side,
in block $(i,j)$ for $1 \le i \le t$ and $1 \le j \le n$, we put $\rho_\lambda (E^i_j)$, the representation
matrix of $E^i_j$.  The matrix $X_\lambda$ is the representation matrix for the components in partition
$\lambda$ of the expansions of the ternary association types in degree $n$. We compute the row canonical form
of $X_\lambda$ and distinguish the upper (respectively lower) part containing the rows with leading ones in
the left (respectively right) side. The rows of the lower right part represent polynomial identities satisfied
by the PATS as a result of dependence relations among the dialgebra expansions of the ternary association
types.

\begin{definition} \label{allrank}
The number of (nonzero) rows in the lower right block of the row canonical form of $X_\lambda$ is the
\emph{rank of identities satisfied by the PATS in partition $\lambda$}.
\end{definition}


\section{Degree 7: nonlinear identities} \label{nonlinear}

In this section we find nonlinear identities for the PATS in degree 7 which are shorter than
the 120-term multilinear identities $R$ and $S$ of Definition \ref{RSdefinition}.

\begin{definition}
A polynomial identity is \emph{nonlinear} if it is homogeneous of degree $n$ and there is a partition
$(n_1, \dots, n_\ell)$ of $n$ with some $n_i \ge 2$ (equivalently $\ell < n$) such that the variables
in each monomial are a permutation of
  \begin{eqnarray*}
  \fl \qquad
  \overbrace{a_1, \dots, a_1}^{n_1}, \,
  \overbrace{a_2, \dots, a_2}^{n_2}, \,
  \dots, \,
  \overbrace{a_\ell, \dots, a_\ell}^{n_\ell}.
  \end{eqnarray*}
The representation theory of the symmetric group tells us which partitions to use in our search for shorter
nonlinear identities.
\end{definition}

\subsection{Application of representation theory}

In Section \ref{degree7section}, we found the inequivalent multilinear monomials in each ternary association
type.  That method was based on monomials and used the skew-symmetries implied by $P \equiv 0$ and $Q \equiv 0$
to reduce the number of monomials in each type.  In contrast, the representation theory is based on association
types and expresses the skew-symmetries as multilinear polynomial identities.

\begin{lemma}
In a partially alternating ternary algebra, every skew-symmetry of the 5 association types in degree 7
is a consequence of the 15 identities in Table \ref{skewsymmetry7}.
\end{lemma}

\begin{proof}
This follows by applying the identities $P \equiv 0$ and $Q \equiv 0$.
\end{proof}

  \begin{table}
  \caption{Skew-symmetries of ternary association types in degree 7.}
  \label{skewsymmetry7}
  \begin{indented}
  \item[]
  $
  \begin{array}{rr}
  \br
  (((abc)de)fg) + (((acb)de)fg) \equiv 0 &\quad (((abc)de)fg) + (((abc)ed)fg) \equiv 0 \\
  (((abc)de)fg) + (((abc)de)gf) \equiv 0 &\quad ((a(bcd)e)fg) + ((a(cbd)e)fg) \equiv 0 \\
  ((a(bcd)e)fg) + ((a(bdc)e)fg) \equiv 0 &\quad ((a(bcd)e)fg) + ((a(bcd)e)gf) \equiv 0 \\
  ((abc)(def)g) + ((acb)(def)g) \equiv 0 &\quad ((abc)(def)g) + ((abc)(edf)g) \equiv 0 \\
  ((abc)(def)g) + ((abc)(dfe)g) \equiv 0 &\quad (a((bcd)ef)g) + (a((cbd)ef)g) \equiv 0 \\
  (a((bcd)ef)g) + (a((bdc)ef)g) \equiv 0 &\quad (a((bcd)ef)g) + (a((bcd)fe)g) \equiv 0 \\
  (a(bcd)(efg)) + (a(cbd)(efg)) \equiv 0 &\quad (a(bcd)(efg)) + (a(bdc)(efg)) \equiv 0 \\
  (a(bcd)(efg)) + (a(efg)(bcd)) \equiv 0 \\
  \br
  \end{array}
  $
  \end{indented}
  \end{table}

  \begin{table}
  \caption{Matrix ranks of each representation in degree 7 for the PATS.}
  \label{degree7ranks}
  \begin{indented}
  \item[]
  \begin{tabular}{rlrrrr}
  \br
  &\quad partition &\quad dimension &\quad symrank &\quad exprank &\quad newrank \\
   1 &\quad 7       &\quad  1 &\quad   5 &\quad   5 &\quad 0 \\
   2 &\quad 61      &\quad  6 &\quad  30 &\quad  30 &\quad 0 \\
   3 &\quad 52      &\quad 14 &\quad  70 &\quad  70 &\quad 0 \\
   4 &\quad 511     &\quad 15 &\quad  75 &\quad  75 &\quad 0 \\
   5 &\quad 43      &\quad 14 &\quad  69 &\quad  69 &\quad 0 \\
   6 &\quad 421     &\quad 35 &\quad 170 &\quad 170 &\quad 0 \\
   7 &\quad 4111    &\quad 20 &\quad  96 &\quad  96 &\quad 0 \\
   8 &\quad 331     &\quad 21 &\quad  99 &\quad  99 &\quad 0 \\
   9 &\quad 322     &\quad 21 &\quad  96 &\quad  96 &\quad 0 \\
  10 &\quad 3211    &\quad 35 &\quad 156 &\quad 156 &\quad 0 \\
  11 &\quad 31111   &\quad 15 &\quad  63 &\quad  64 &\quad 1 \\
  12 &\quad 2221    &\quad 14 &\quad  56 &\quad  56 &\quad 0 \\
  13 &\quad 22111   &\quad 14 &\quad  52 &\quad  53 &\quad 1 \\
  14 &\quad 211111  &\quad  6 &\quad  17 &\quad  20 &\quad 3 \\
  15 &\quad 1111111 &\quad  1 &\quad   0 &\quad   2 &\quad 2 \\
  \br
  \end{tabular}
  \end{indented}
  \end{table}

Let $\lambda$ be a partition of $n = 7$ with associated irreducible representation of dimension
$d = d_\lambda$. There are 15 skew-symmetry identities and 5 association types, requiring a
$15d \times 5d$ matrix $M_\lambda$. In each skew-symmetry, the first term has the identity permutation
with representation matrix $I_d$, and the second term has a permutation $\pi$ of order 2 with
representation matrix $(R^\lambda_\mathrm{id})^{-1} R^\lambda_\pi$ from Lemma \ref{cliftonlemma}.
The $d \times d$ block in position $(i,j)$ contains the sum of these matrices, where $i$ and $j$ are
the index number and the association type of the skew-symmetry. The rank of $M_\lambda$ is `symrank'
in Table \ref{degree7ranks}.  For each $\lambda$, we construct the representation matrix $X_\lambda$
(Table \ref{repmatexp}) and compute the rank of its lower right part (Definition \ref{allrank}); this
is `exprank' in Table \ref{degree7ranks}.  Column `newrank' is the difference between `symrank' and
`exprank': this is the rank of the new identities in degree 7 for partition $\lambda$; that is, the
identities which are not trivial consequences of the skew-symmetries of the association types. We check
the results by summing, over all representations, the product of `newrank' and `dimension':
$1 \cdot 15 + 1 \cdot 14 + 3 \cdot 6 + 2 \cdot 1 = 49$. This is the dimension of the nullspace of the
expansion matrix from Section \ref{degree7section}; column `newrank' gives the decomposition of the
nullspace into irreducible representations.

There are four representations where `newrank' is positive: 11, 13, 14, 15. This suggests that a slight
modification of the techniques of Section \ref{degree7section} will produce nonlinear identities in these
partitions: identities in which the variables in each term are a permutation of $\{a, a, a, b, c, d, e\}$,
$\{a, a, b, b, c, d, e\}$ or $\{a, a, b, c, d, e, f\}$. (We omit representation 15 since it corresponds to
the multilinear case.)

\subsection{Straightening algorithm} \label{subsectionstraightening}

We need an algorithm to convert a monomial (multilinear or nonlinear) to its `straightened' form with respect
to $P \equiv 0$ and $Q \equiv 0$. We apply $P$ and $Q$ to convert a monomial with a given permutation of the
variables into ($\pm$) a monomial with a different permutation which lexicographically precedes the original
permutation.
We use the recursiveprocedures \texttt{completestraighten} (\texttt{CS}) and \texttt{partialstraighten}
(\texttt{PS}); for both the input is a monomial $x$. If the straightened form of $x$ is 0, then both return 0.
If $\deg(x) = 1$ then both return $x$. If $\deg(x) > 1$ then we write $x = ( x_1, x_2, x_3 )$ and proceed as
follows:
  \begin{itemize}
  \item
  \texttt{CS} recursively computes $\texttt{CS}(x_1)$, $\texttt{CS}(x_2)$, $\texttt{CS}(x_3)$. If any of
  them is 0, then \texttt{CS} returns 0. If two or more are equal, then \texttt{CS} returns 0. Otherwise,
  \texttt{CS} puts them in the correct order using \texttt{strictlyprecedes} (see below).
  \item
  \texttt{PS} recursively computes $\texttt{PS}(x_1)$, $\texttt{CS}(x_2)$, $\texttt{CS}(x_3)$.  If any of
  them is 0, then \texttt{PS} returns 0. If $\texttt{CS}(x_2) = \texttt{CS}(x_3)$, then \texttt{PS} returns
  0. Otherwise, \texttt{PS} puts $\texttt{CS}(x_2)$ and $\texttt{CS}(x_3)$ in the correct order using
  \texttt{strictlyprecedes}.
  \end{itemize}
Procedure \texttt{strictlyprecedes} compares monomials $x$ and $y$. If $\deg(x) \ne \deg(y)$ then it
returns \texttt{true} if $\deg(x) < \deg(y)$, \texttt{false} if $\deg(x) > \deg(y)$.
If $\deg(x) = \deg(y)$ then:
  \begin{itemize}
  \item
  If both $x$ and $y$ have degree 1, it uses the total order on the generators.
  \item
  If both have degree $> 1$ then $x = ( x_1, x_2, x_3 )$ and $y = ( y_1, y_2, y_3 )$ and it finds the least
  $i$ with $x_i \ne y_i$ and recursively calls \texttt{strictlyprecedes}$(x_i,y_i)$.
  \end{itemize}
In other words, first compare the degrees;
if the degrees are equal then compare the association types;
and if the types are equal then compare the permutations.

\begin{theorem} \label{nonlinear31111}
There is one identity for partition 31111; it has 60 terms:
  \begin{eqnarray*}
  \fl \qquad
  I_{31111}^{(1)}
  &=
  \frac14
  \sum_{\sigma \in S_5}
  \epsilon(\sigma) \,
  [ [ [ a a^\sigma b^\sigma ] a c^\sigma] d^\sigma e^\sigma ]
  -
  \frac1{12}
  \sum_{\sigma \in S_5}
  \epsilon(\sigma) \,
  [ [ a [ a^\sigma b^\sigma c^\sigma ] a ] d^\sigma e^\sigma ]
  \\
  \fl \qquad
  &\quad
  -
  \frac1{12}
  \sum_{\sigma \in S_5}
  \epsilon(\sigma) \,
  [ [ a a^\sigma b^\sigma ] [ c^\sigma d^\sigma e^\sigma ] a ]
  -
  \frac16
  \sum_{\sigma \in S_5}
  \epsilon(\sigma) \,
  [ a [ [ b^\sigma c^\sigma d^\sigma ] a e^\sigma ] a^\sigma ].
  \end{eqnarray*}
\end{theorem}

\begin{proof}
We first generate all 840 permutations of $a, a, a, b, c, d, e$.  For each association type, we apply the type
to each permutation, find the straightened form of the resulting monomial, and retain only those monomials
which equal their own straightened forms.  We sort the remaining monomials in each type by lex order of the
permutation.  For partition 31111, the five association types contain respectively $60 + 34 + 34 + 34 + 3 =
165$ monomials.  The expansion matrix $E$ has 165 columns and $7 \cdot 840 = 5880$ rows.
For $j = 1, \dots, 165$ we store in column $j$ the PATS expansion of ternary monomial $j$.  The rank is 164,
and so the nullspace has dimension 1.  Hence, up to a scalar multiple, there is exactly one identity for the
PATS with these variables; this identity has 60 terms with coefficients $\pm 1$ in association types 1--4.
\end{proof}

\begin{theorem} \label{nonlinear22111}
There are two identities for partition 22111; both have 60 terms:
  \begin{eqnarray*}
  \fl \qquad
  I_{22111}^{(1)}
  &=
  \frac14
  \sum_{\sigma \in S_5}
  \epsilon(\sigma) \,
  [ [ [ a a^\sigma b^\sigma] b c^\sigma] d^\sigma e^\sigma ]
  -
  \frac{1}{12}
  \sum_{\sigma \in S_5}
  \epsilon(\sigma) \,
  [ [ a (a^\sigma b^\sigma c^\sigma] b] d^\sigma e^\sigma]
  \\
  \fl \qquad
  &\quad
  -
  \frac{1}{12}
  \sum_{\sigma \in S_5}
  \epsilon(\sigma) \,
  [ [ a a^\sigma b^\sigma] [c^\sigma d^\sigma e^\sigma] b]
  -
  \frac16
  \sum_{\sigma \in S_5}
  \epsilon(\sigma) \,
  [ a [[ a^\sigma b^\sigma c^\sigma] b d^\sigma] e^\sigma],
  \end{eqnarray*}
and $I_{22111}^{(2)}$ which is obtained by interchanging $a$ and $b$.
\end{theorem}

\begin{proof}
Similar to the proof of Theorem \ref{nonlinear31111}.
\end{proof}

\begin{theorem} \label{nonlinear211111}
There are 12 identities for partition 211111; only 5 have 60 terms:
  \begin{eqnarray*}
  \fl \qquad
  I_{211111}^{(1)}
  &
  =
  \,
  \frac14
  \sum_{\sigma \in S_5}
  \epsilon(\sigma) \,
  [ [ [ b a^\sigma c^\sigma] a d^\sigma] e^\sigma f^\sigma]
  -
  \frac{1}{12}
  \sum_{\sigma \in S_5}
  \epsilon(\sigma) \,
  [ [ b [ a^\sigma c^\sigma d^\sigma] a ] e^\sigma f^\sigma]
  \\
  \fl \qquad
  &\quad
  -
  \frac{1}{12}
  \sum_{\sigma \in S_5}
  \epsilon(\sigma) \,
  [ [ b a^\sigma c^\sigma] [d^\sigma e^\sigma f^\sigma] a]
  -
  \frac16
  \sum_{\sigma \in S_5}
  \epsilon(\sigma) \,
  [ b [[ a^\sigma c^\sigma d^\sigma] a e^\sigma] f^\sigma],
  \end{eqnarray*}
and $I_{211111}^{(2)}, \dots, I_{211111}^{(5)}$ obtained by interchanging $b$ and $c$, $d$, $e$, $f$
respectively.
\end{theorem}

\begin{proof}
Similar to the proof of Theorem \ref{nonlinear31111}.
\end{proof}


\section{Degree 9} \label{degree9section}

In degree 9 there are 12 association types for a PATA (Table \ref{CAPAtypes}).
We compute the following matrix ranks using the representation theory of the symmetric group:
  \begin{itemize}
  \item symrank:
  the rank of the skew-symmetry identities of the association types.
  \item symlifrank:
  the rank of the skew-symmetries combined with the consequences in degree 9 of the identities $R$ and $S$
  in degree 7 from Definition \ref{RSdefinition}; that is,
  \begin{eqnarray*}
  \fl \qquad
  T((ahi),b,c,d,e,f,g), \;
  T(a,(bhi),c,d,e,f,g), \;
  \dots, \;
  T(a,b,c,d,e,f,(ghi)),
  \\
  \fl \qquad
  (T(a,b,c,d,e,f,g),h,i), \;
  (h,T(a,b,c,d,e,f,g),i), \;
  (h,i,T(a,b,c,d,e,f,g)),
  \end{eqnarray*}
  where $T = R$ and $T = S$.
  \item exprank:
  the rank of the lower right part of the expansion matrix (Definition \ref{allrank}).
  \end{itemize}
For every partition `symlifrank' equals `exprank': there are no new identities.


\section{Conclusion}

Trilinear operations in an associative algebra have recently been classified by Bremner and Peresi
\cite{BP-CTO}: there are six isolated operations (the alternating, symmetric, and cyclic sums, the cyclic
commutator, the weakly commutative and anticommutative operations), and four infinite families (the Lie,
Jordan and anti-Jordan families, and a fourth family which seems unrelated to Lie and Jordan structures).
The Kolesnikov-Pozhidaev algorithm can be applied to all these operations:
we choose one of the three arguments and make it the center of each monomial.
Unlike the alternating sum, which corresponds to the 1-dimensional sign representation of the
symmetric group $S_3$, most of the other algebra operations will produce essentially
different dialgebra operations for each choice of the central argument.
The polynomial identities satisfied by these new dialgebra operations will define varieties of ternary
algebras with great potential for applications in pure mathematics and theoretical physics.


\ack
The first author was supported by a Discovery Grant from NSERC,
the Natural Sciences and Engineering Research Council of Canada.
The second author was supported by the Spanish MEC and Fondos FEDER jointly
through project MTM2007-60333, and by the Junta de Andaluc\'ia (projects FQM-336 and FQM2467).
She thanks the Department of Mathematics and Statistics at the University of Sas\-kat\-che\-wan
for its hospitality in June 2010, and the President's Diversity Enhancement Fund
and the Role Model Speaker Fund of the College of Arts and Science for financial support.
Both authors thank Hader Elgendy for the reference to the survey \cite{AI}.


\end{document}